\documentclass[a4paper,11pt]{article}

\usepackage{amssymb,amsmath,amsthm}
\usepackage{epsfig}
\usepackage{breqn}
\usepackage{rotating}
\usepackage{amsfonts}
\usepackage{setspace}
\usepackage{fullpage}
\usepackage{enumitem}
\usepackage{bbold} 
\usepackage{amsmath}
\usepackage{comment}
\usepackage{pgf,tikz}
\usepackage{mathtools}  
\usepackage{hyperref}

\usepackage{authblk}

\usepackage{setspace}

\newtheorem{theorem}{Theorem}
\newtheorem{claim}[theorem]{Claim}

\newtheorem{definition}[theorem]{Definition}

\newcommand{\floor}[1]{\left\lfloor{#1}\right\rfloor}
\newcommand{\ceil}[1]{\left\lceil{#1}\right\rceil}

\def\C{\mathcal C}
\def\B{\mathcal B}
\def\h{\mathcal H}

\def\F{\mathcal F}

\newcommand{\abs}[1]{\left\lvert{#1}\right\rvert}

\def\path{semi-path}

\begin{document}

\title{Connected Hypergraphs without Long Berge-Paths}

\author[1,2]{Ervin Gy\H{o}ri} 

\author[1,2]{Nika Salia}
\author[2,4]{Oscar Zamora} 

\affil[1]{Alfr\'ed R\'enyi Institute of Mathematics, Budapest. \par 
 \texttt{gyori.ervin@renyi.hu, salia.nika@renyi.hu}}
\affil[2]{Central European University, Budapest.     }
 \affil[4]{Universidad de Costa Rica, San Jos\'e. \par
 \texttt{oscar.zamoraluna@ucr.ac.cr}}

\maketitle

\begin{abstract}

We generalize a result of Balister, Gy{\H{o}}ri, Lehel and Schelp for hypergraphs. We determine the unique extremal structure  of an $n$-vertex $r$-uniform connected  hypergraph with the maximum number of hyperedges, without a $k$-Berge-path, where $n \geq N_{k,r}$, $k\geq 2r+13>17$.

\end{abstract}

\section{Basic Definitions with Background.}

Tur\'an-type problems are central in extremal combinatorics. Classic Tur\'an questions are to determine  the maximum number of edges in an $n$ vertex graph without certain substructure. We give a formal definition of the Tur\'an number.
 
\begin{definition}
The Tur\'an number of a graph $H$, denoted by $ex(n,H)$, is the maximum number of edges in a graph on n vertices that do not contain $H$, as a subgraph. The
Tur\'an number of a family of graphs $\h$, denoted by $ex(n,\h)$, is the maximum number of edges in a graph on n vertices that do not contain any $H$, $H \in \h$, as a subgraph.
\end{definition}

Let us  denote by  $P_k$ the  path of length $k$ and we denote by $\C_{\geq k}$ the family of cycles of length at least $k$.  Erd\H{o}s and Gallai \cite{Er-Ga} determined the Tur\'an number  of paths and cycles. Here we state those results. 

\begin{theorem}[Erd\H{o}s, Gallai~\cite{Er-Ga}]\label{EG1} 
For all integers $n$ and $k$ satisfying $n \geq k \geq 1$, we have   
\[ex(n,P_k)\leq \frac{(k - 1)n}{2}.\]
\end{theorem}

Erd\H{o}s and Gallai \cite{Er-Ga} proved the following result on the Tur\'an number of cycles.

\begin{theorem}[Erd\H{o}s, Gallai~\cite{Er-Ga}]\label{EG2} 
For all integers $n$ and $k$ satisfying $n \geq k \geq 3$,  we have 
\[ex(n,\C_{\geq k}) \leq \frac{(k - 1)(n - 1)}{2}.\]
\end{theorem}

In fact Theorem \ref{EG1} was deduced as a simple corollary of Theorem \ref{EG2}. Later, in 1975, this result was improved by Faudree and Schelp \cite{FS}. They have determined  $ex(n,P_k)$ for all $n$ and $k$, satisfying $n>k>0$, with
the corresponding extremal graphs. 

In 1977, Kopylov, \cite{Kopylov},  determined $ex^{conn}(n,P_k)$, where  $ex^{conn}(n,\cdot)$ denotes classical Tur\'an  number for connected graphs. Surprisingly the asymptotics was the same as in the non-connected case for all paths of odd length. In 2008, Balister, Gy{\H{o}}ri, Lehel and Schelp \cite{BGLS} improved Kopylov's result by characterizing the extremal graphs for all $n$.  Recently numerous mathematicians started investigating similar problems for $r$-uniform hypergraphs.  Before presenting some of those results, related to our work, we need to state some useful definitions. In what follows, we refer to simple $r$-uniform hypergraphs as $r$-graphs.


\begin{definition} [Berge \cite{Berge}]
A \emph{Berge-path} of length $t$  is an alternating sequence of distinct $t+1$ vertices and distinct $t$ hyperedges of the hypergraph,  $v_1, e_1, v_2, e_2, v_3, \dots, e_t, v_{t+1}$ such that, $v_{i},v_{i+1} \in e_i$, for $i =1,2,\dots, t$. We denote the set of all Berge-paths of length $t$ by $\B P_t$. Vertices $v_1, v_2, \dots, v_{t+1}$ are called defining vertices and hyperedges $e_1,e_2,\dots,e_t$ are called defining hyperedges of the path.  
\end{definition}

\begin{definition} [Berge \cite{Berge}]
A \emph{Berge-cycle} of length $t$ is an alternating sequence of $t$ distinct vertices and $t$ distinct hyperedges of the hypergraph,  $v_1, e_1, v_2, e_2, v_3, \dots, v_t, e_t$, such that $v_{i},v_{i+1} \in e_i$, for $i =1,2,\dots, t$, indices are taken modulo $t$.  Let us denote set of all Berge-cycles of length at least $t$ by $\B C_{\geq t}$. Vertices $v_1, v_2, \dots, v_{t}$ are called defining vertices and hyperedges $e_1,e_2,\dots,e_t$ are called defining hyperedges of the cycle.   
\end{definition}
 


It is natural to ask Tur\'an-type questions for the hypergraphs.  The Tur\'an number of a family of $r$-graphs $\h$, denoted by $ex_r(n,\h)$, is the maximum number of hyperedges in an $r$-graph on n vertices that does not contain any copy of $H$ as a sub-hypergraph, for all $H \in \h$. Moreover one may require connectivity and ask  Tur\'an-type questions for connected $r$-graphs. We denote by $ex^{conn}_r(n,\h)$ the maximum number of hyperedges in a connected $r$-graph on n vertices that does not contain any copy of $H$ as a sub-hypergraph for all $H \in \h$.




The first extension of Erd\H{o}s and Gallai theorems \cite{Er-Ga} for $r$-graphs was by Gy\H{o}ri, Katona, and Lemons~\cite{GyoKaLe}, in 2016 they extended Theorem \ref{EG1} to $r$-graphs. It turned out that the extremal hypergraphs have a different behavior when $k \leq r$ and $k > r$.
\begin{theorem}[Gy\H{o}ri, Katona and Lemons~\cite{GyoKaLe}]\label{GKL1} For all integers $n$, $r$, and $k$ satisfying $r \geq k \geq 3$, we have 
\[ex_r(n,\B P_{k}) \leq \frac{(k-1)n}{r+1}.\]
The equality holds if and only if $(r+1)|n$, and the extremal $r$-graph is the disjoint union of $\frac{n}{r+1}$  sets of size $r+1$ containing $k-1$ hyperedges each.
\end{theorem}
\begin{theorem}[Gy\H{o}ri, Katona and Lemons~\cite{GyoKaLe}]\label{GKL2}  For all integers $n$, $r$, and $k$ satisfying $k > r + 1 > 3$, we have 
\[ex_r(n,\B P_{k}) \leq \frac{n}{k}\binom{k}{r}.\] 
The equality holds if and only if $k|n$, and the only extremal $r$-graph is the disjoint union of $\frac{n}{k}$ copies of the complete $k$-vertex $r$-graph.
\end{theorem}
The remaining case when $k = r + 1$ was solved later by Davoodi, Gy\H{o}ri, Methuku, and Tompkins~\cite{DavoodiGMT}. They showed that the extremal number matches the upper bound of Theorem  \ref{GKL2}.

Observe that the extremal hypergraph is not connected in Theorem \ref{GKL1} and Theorem \ref{GKL2}. Naturally, one may require connectivity and search for connected extremal hypergraphs. The  extension of Theorem \ref{EG2}, for Berge $r$-graphs, was solved by  different groups of researchers recently. As expected the extremal hypergraphs were different for the cases $r+1<k$ and $r+1 \geq k$. Moreover the case $k=r$ presents a different behavior. The investigation of the Tur\'an number of the long Berge-cycles started with the paper of F\"uredi, Kostochka and Luo \cite{furedi2018avoiding,furedi2018avoiding2}. They got the sharp result for all $n$, when $k \geq r + 3$. The case when $k=r + 2$ together with the case $k=r+1$ was settled later by authors of this paper together with  Ergemlidze, Methuku and Tompkins ~\cite{ergemlidze2018avoiding}. Those results were followed by  Kostochka and Luo \cite{KostochkaLuo}, for $r \geq k + 1$,  which left out the last  mysterious case where $k$ and $r$ are equal. The missing case  $k=r$, together with  $r \geq k + 1$ was settled by authors of this paper together with Lemons  \cite{GLSZ}.

\begin{theorem}[F\"uredi, Kostochka and Luo~\cite{furedi2018avoiding,furedi2018avoiding2}] \label{t1} For all integers $n$, $r$, and $k$ satisfying  $k \geq r + 3 \geq 6$, we have

\[ex_r(n,\B C_{\geq k})\leq \frac{n-1}{k-2}\binom{k-1}{r}.\]

The equality holds if and only if $(k-2)|(n-1)$ and the extremal $r$-graph are the  union of $\frac{n-1}{k-2}$ copies of complete $(k-1)$-vertex $r$-graphs sharing a vertex in a  tree-like structure.
\end{theorem}

The first attempt, to determine the  Tur\'an number, for connected $r$-graphs without a Berge-path of length $k$ was published in 2018,  \cite{GyoriMSTV}. They  determined the asymptotics of the extremal function. 

Gy{\H{o}}ri, Methuku, Salia, Tompkins and Vizer determined asymptotics of the  Tur\'an number for connected $r$-graphs without a Berge-path of length $k$.

\begin{theorem}[Gy{\H{o}}ri, Methuku, Salia, Tompkins, Vizer \cite{GyoriMSTV}]\label{Thm:con_ass}
Let $\h_{n,k}$ be a largest,  $r$-uniform, connected, $n$-vertex hypergraph with no Berge-paths of length $k$, then
\begin{displaymath}
\lim_{k \to \infty} \lim_{n \to \infty}  \frac{\abs{E(\h_{n,k})}}{k^{r-1}n} = \frac{1}{2^{r-1}(r-1)!}.
\end{displaymath} 
\end{theorem}

In the recent work F\"uredi, Kostochka and Luo  investigated the $2$-connected hypergraphs and obtained a number of interesting extensions of Theorems \ref{EG1} and Theorem \ref{EG2}. Before stating the result strengthening Theorem \ref{Thm:con_ass}, we define the following function $\displaystyle f^*(n,k,r,a):={\binom{k-a}{r}} + (n-k+a) {\binom{a} {r-1}}$.

\begin{theorem}[F\"uredi, Kostochka, Luo, ~\cite{ConnFKL}] \label{FKL_EG_Conn}
Let $n \geq n'_{k,r} \geq k \geq 4r\geq 12$. If $\h$ is an $n$-vertex connected $r$-graph with no Berge-paths of length $k$, then $$e(\h) \leq  f^*(n,k,r, \lfloor (k-1)/2 \rfloor).$$ 
\end{theorem}

In the following section we  prove Theorem \ref{FKL_EG_Conn} for $k \geq 2r+13$. We also prove that  there is a unique hypergraph with the extremal number of hyperedges. The rough idea of the proof is to investigate the structure of an $r$-graph with the extremal properties.  We provide the key definition used in the following section. The vertex neighborhood of a vertex $v$ in an $r$-graph $\h$, is denoted by $N_{\h}(v)$, $N_{\h}(v):=\Big\{u \in V(\h)\Big| \exists h, h \in \h \mbox{ and } \{v,u\}\in h \Big\}$.


\begin{figure}
\label{figureConstruction}
\centering
\begin{tikzpicture}
\draw (0,0) arc(0:360:1cm and .5cm) (1.1,-2) arc(0:360:2.1cm and .9cm);

\foreach \x in {-2.75,-2.05,...,0.75}{
\filldraw[blue, fill opacity=0.10] plot [smooth cycle] coordinates {(-1.8,0) (-1.2,.3) (-.6,0)(\x+0.2,-2.1) (\x-0.2,-2.1)  } ;
}



\foreach \x in {-2.75,-2.05,...,0.75}{
\filldraw[blue, fill opacity=0.10] plot [smooth cycle] coordinates {(-1.4,0) (-.8,.3) (-.2,0)(\x+0.2,-2.1) (\x-0.2,-2.1)  } ;
\filldraw[red] (\x,-2) circle (2pt);
}

\filldraw[red] (-1.6,0) circle (2pt) (-1,0) circle (2pt) (-.4,0) circle (2pt);

\node at (-1,-3.15){If $k$ is odd.};
\end{tikzpicture}\qquad
\begin{tikzpicture}
\draw (0,0) arc(0:360:1cm and .5cm) (-.2,-2) arc(0:360:1.5cm and .7cm);

\foreach \x in {-2.75,-2.05,...,-0.65}{
\filldraw[blue, fill opacity=0.10] plot [smooth cycle] coordinates {(-1.8,0) (-1.2,.3) (-.6,0)(\x+0.2,-2.1) (\x-0.2,-2.1)  } ;
}

\foreach \x in {0.5,1.25}{
\filldraw[blue, fill opacity=0.10] plot [smooth cycle] coordinates {(-1.8,0) (-1.2,.3) (-.6,0)(\x+0.2,-2.1) (\x-0.2,-2.1)  } ;
}



\foreach \x in {-2.75,-2.05,...,-0.65}{
\filldraw[blue, fill opacity=0.10] plot [smooth cycle] coordinates {(-1.4,0) (-.8,.3) (-.2,0)(\x+0.2,-2.1) (\x-0.2,-2.1)  } ;
\filldraw[red] (\x,-2) circle (2pt);
}

\foreach \x in {0.5,1.25}{
\filldraw[blue, fill opacity=0.10] plot [smooth cycle] coordinates {(-1.4,0) (-.8,.3) (-.2,0)(\x+0.2,-2.1) (\x-0.2,-2.1)  } ;
\filldraw[red] (\x,-2) circle (2pt);
}

\filldraw[red] (-1.6,0) circle (2pt) (-1,0) circle (2pt) (-.4,0) circle (2pt);

\filldraw[blue, fill opacity=0.10] plot [smooth cycle] coordinates {(0.9+.5,-2) (.3+.5,-2.2) (-0.2+.5,-2) (-1.6-0.2,0.1) ( (-1.6+0.2,0.1)  } ;

\filldraw[blue, fill opacity=0.10] plot [smooth cycle] coordinates {(0.9+.5,-2) (.3+.5,-2.2) (-0.2+.5,-2) (-1.2-0.1,0.1) (-1.2+0.3,0.1) } ;

\filldraw[blue, fill opacity=0.10] plot [smooth cycle] coordinates {(0.9+.5,-2) (.3+.5,-2.2) (-0.2+.5,-2) (-.6-0.1,0.1) (-.6+0.3,0.1) } ;

\node at (-1,-3.15){If $k$ is even.};
\end{tikzpicture}
\caption{Graphs $\h_{n,k}$}
\end{figure}


\section[]{Proof of the uniqueness of  the extremal, connected, $r$-graph, without a long Berge-path.}

Let us start this section by defining $\h_{n,k}$, the extremal hypergraph for the following theorem. There are two similar $n$ vertex, $r$-graphs, without a Berge-path of length $k$, depending on the parity of $k$.   If $k$ is odd, the vertex set is partitioned into two sets $A$ and $B$, $\abs{A}=\frac{k-1}{2}$ and $\abs{B}=n-\abs{A}$.  The hyperedges are all $r$-sets containing at most one vertex from the set $B$ the rest from the set $A$. If $k$ is even the vertex set is partitioned into two sets $A$ and $B$, with sizes $\abs{A}=\frac{k-2}{2}$ and $\abs{B}=n-\abs{A}$. Two distinct vertices are fixed in the set $B$, say $b_1$ and $b_2$. The hyperedges are all $r$-sets containing at most one vertex from the set $B$ or containing $b_1$ and $b_2$ and the rest from $A$ see Figure~\ref{figureConstruction}.  Please note that in order to have a hyperedge in this graph it is necessary to have $\floor{\frac{k-1}{2}}\geq{r-1}$.

\begin{theorem}\label{main_theorem}
For all integers  $n$, $k$, and $r$  there exists $N_{k,r}$ such that if $n>N_{k,r}$ and $k \geq 2r+13 \geq 18$, then
 
 \begin{displaymath}
 ex_r^{conn}(n,\B P_k)={\binom{\floor{\frac{k-1}{2}}}{r-1}}\bigg(n-\floor{\frac{k-1}{2}}\bigg)+{\binom{\floor{ \frac{k-1}{2}}}{r}}+ \mathbb{1}_{2|k}{\binom{\floor{\frac{k-1}{2}}}{r-2}}.
 \end{displaymath}
 
The extremal hypergraph is the unique $r$-graph, $\h_{n,k}$.
\end{theorem}

For the existence of the $r$-graph, $\h_{n,k}$, we need $\floor{\frac{k-1}{2}}\geq{r-1}$, therefore the condition $k \geq 2r+13$ is natural and asymptotically best possible.  We use this lower bound for the case when $r=3$. It can be improved slightly for larger $r$. 

\begin{proof}
The lower bound comes from $\h_{n,k}$. Observe that from  two consecutive defining vertices, of a Berge-path in $\h_{n,k}$,  at least one is from the set $A$ or they are  $b_1$ and $b_2$. Therefore  the length of a longest Berge-path of $\h_{n,k}$ is at most $k-1$. 
 An $r$-graph $\h_{n,k}$, is connected and contains no Berge-path of length $k$. Therefore, we have a desired lower bound  $\displaystyle {\binom{\floor{\frac{k-1}{2}}}{r-1}}\bigg(n-\floor{\frac{k-1}{2}}\bigg)+{\binom{\floor{ \frac{k-1}{2}}}{r}}+ \mathbb{1}_{2|k}{\binom{\floor{\frac{k-1}{2}}}{r-2}}$. 

For the upper bound and uniqueness of the extremal $r$-graph, we consider an  $r$-graph on $n$ vertices, without a Berge-path of length $k$, with at least $\displaystyle{\binom{\floor{\frac{k-1}{2}}}{r-1}}\left(n-\floor{\frac{k-1}{2}}\right)+{\binom{\floor{ \frac{k-1}{2}}}{r}}+ \mathbb{1}_{2|k}{\binom{\floor{\frac{k-1}{2}}}{r-2}}$ hyperedges, denoted by $\h$. 
For large enough $n$, namely  if $\displaystyle n> N_{k,r}:=\binom{r(k-1)}{r}+{\binom{\floor{\frac{k-1}{2}}}{r-1}} \floor{\frac{k-1}{2}}-{\binom{\floor{ \frac{k-1}{2}}}{r}}- \mathbb{1}_{2|k}{\binom{\floor{\frac{k-1}{2}}}{r-2}} +r \cdot (k-1),$ there exists an $r$-uniform, sub-hypergraph $\h'$ of $\h$ with  at least $r\cdot (k-1)$ vertices and the minimum degree at least 
$\displaystyle{\binom{\floor{\frac{k-1}{2}}}{r-1}}$. The $r$-graph $\h'$ is obtained from $\h$ by repeatedly removing vertices of  degree less than $\displaystyle{\binom{\floor{\frac{k-1}{2}}}{r-1}}$ and  all hyperedges incident to them.

\begin{claim}\label{size_of_shadow}
There exists a longest Berge-path $P$ in $\h'$, with  defining vertices $U=\{u_1, \dots, u_{l+1}\}$ and defining hyperegdes $\F=\{f_1,f_2, \dots, f_{l}\}$ in this given order, such that the sizes of the vertex neighborhood of $u_1$ and $u_{l+1}$ in $\h' \setminus \F$ are at least $\floor{\frac{k-1}{2}}$, $\abs{N_{\h' \setminus \F}(u_1)}\geq \floor{\frac{k-1}{2}}$ and $\abs{N_{\h' \setminus \F}(u_{l+1})}\geq \floor{\frac{k-1}{2}}$.
\end{claim}
\begin{proof}[Proof of Claim \ref{size_of_shadow}]
Let $P$ be a longest Berge-path in $\h'$, minimizing $x_1+x_{l+1}$, where  $x_1$ and $x_{l+1}$ are the number of hyperedges from $\F$ incident with $u_1$ and $u_{l+1}$ respectively. If $u_1 \in f_i$ then the Berge-path 
\begin{displaymath}
u_i, f_{i-1},u_{i-1},f_{i-2}, u_{i-2}, \dots, u_2, f_1, u_1, f_i,u_{i+1}, f_{i+1}, \dots, u_l, f_l, u_{l+1}
\end{displaymath}
is  a longest Berge-path, with the same set of defining vertices and defining hyperedges. Hence by the minimality of the sum $x_1+x_{l+1}$, the number of  hyperedges from $\F$ incident to $u_i$ is at least $x_1$.

If we consider all possible Berge-paths obtained from $P$ by the way described above (including itself), then the number of pairs $(u,f)$, where $u \in U$, $f \in \F$ and $u \in f$, is at least $x_1^2$, this number is upper bounded by  $r\abs{\F}\leq rl$, since each hyperedge in $\F$ is incident with $r$ vertices.  Hence we have  $x_1^2 \leq rl \leq r(k-1)$, therefore $x_1\leq \sqrt{r(k-1)}$, the same holds for the  other terminal vertex $x_{l+1}$. 
Since the degree of $u_1$ is at least $\displaystyle{\binom{\floor{\frac{k-1}{2}}}{r-1}}$, out of which at most $\sqrt{(k-1)r}$ of the hyperedges are defining hyperedges, we have the degree of $u_1$ in $\h' \setminus \F$ is at least $\displaystyle{\binom{\floor{\frac{k-1}{2}}}{r-1}} - \sqrt{(k-1)r}
> {\binom{\floor{\frac{k-1}{2}}-1}{r-1}}$ this inequality holds, since $k \geq 2\cdot r+13.$ 
Hence we have the desired inequalities $\abs{N_{\h' \setminus \F}(u_1)}\geq \floor{\frac{k-1}{2}}$ and $\abs{N_{\h' \setminus \F}(u_{l+1})}\geq \floor{\frac{k-1}{2}}$.

Please, note that we may use this claim for every hypergraph with same minimum degree condition.
\end{proof}


After removing low degree vertices from $\h$, we may have destroyed one of the important conditions, connectivity. Therefore in the next claim we are going to show that we did not loose this condition. 
\begin{claim}\label{connected}
 $H'$ is connected.
\end{claim}
\begin{proof}[Proof of Claim \ref{connected}]
Assume by contradiction that $\h'$ is not connected. Take a longest Berge-path in two connected components and apply Claim \ref{size_of_shadow} for each of the components.
The size of the vertex neighborhood  without using defining hyperedge of the terminal  vertices of those paths are at least $\floor{\frac{k-1}{2}}$, all of those vertices must be defining vertices, otherwise, we could have extended the longest Berge-path in the components.
Therefore there is a Berge-cycle of length at least   $\floor{\frac{k-1}{2}}+1$ in each of the components. Those Berge-cycles are also in the connected hypergraph  $\h$ since $\h'$ is sub-hypergraph of $\h$.

If $k$ is odd, there is a  Berge-path of length  at least $k$ in $\h$,  containing all vertices of two Berge-cycles from two different components of $\h'$ of length at least $\frac{k-1}{2}+1$ each, a contradiction. 
Similarly, we get a contradiction if $k$ is even and the length of a Berge-cycle in one of the components is at least $\frac{k}{2}+1$.
If $k$ is even and $\h'$ does not contain a Berge-cycle of length at least $\frac{k}{2}+1$, then the number of hyperedges in $\h'$ is at most $\displaystyle\frac{v(\h')-1}{\frac{k}{2}-1} {\binom{\frac{k}{2}}{r}}$ from  Theorem~\ref{t1}.
Since, we  obtained $\h'$ from $\h$ by deleting less than $\displaystyle{\binom{\frac{k}{2}-1}{r-1}}$ hyperedges on each step, we have 

\begin{displaymath}
e(\h) \leq e(\h')+\left({\binom{\frac{k}{2}-1}{r-1}}-1\right)(n-v(\h'))<
\frac{v(\h')-1}{\frac{k}{2}-1} {\binom{\frac{k}{2}}{r}}+{\binom{\frac{k}{2}-1}{r-1}}(n-v(\h')) 
\end{displaymath}
\begin{displaymath}
=\frac{v(\h')-1}{\frac{k}{2}-1}\cdot\frac{\frac{k}{2}}{r}\cdot{\binom{\frac{k}{2}-1}{r-1}}+{\binom{\frac{k}{2}-1}{r-1}}(n-v(\h'))
\end{displaymath}
\begin{displaymath}
<{\binom{\frac{k}{2}-1}{r-1}}\bigg(n-\frac{k-2}{2}\bigg)+{ \binom{\frac{k-2}{2}}{r}}+ {\binom{\frac{k-2}{2}}{r-2}}=e(\h_{n,k}).
\end{displaymath}
We got  $e(\h)< e(\h_{n,k})$ a contradiction, therefore $\h'$ is connected.
\end{proof}


\begin{claim}\label{Contains_cycle}
The hypergraph $\h'$ contains a Berge-cycle of length $l$ or $l+1$ (recall, $l$ is the length of longest Berge-path in $\h'$).
\end{claim}
\begin{proof}[Proof of Claim \ref{Contains_cycle}]

Let $P$ be a longest Berge-path from Claim \ref{size_of_shadow}. Hence we have $\displaystyle\abs{N_{\h' \setminus \F}(u_1)}\geq \floor{\frac{k-1}{2}}$ and $\displaystyle\abs{N_{\h' \setminus \F}(u_{l+1})}\geq \floor{\frac{k-1}{2}}$, where  $u_1,f_1,u_2,\dots,u_{l},f_l,u_{l+1}$ is a longest  Berge-path in $\h'$.
For simplicity, let us denote $S_{u_1}:=N_{\h' \setminus \F}(u_1)$ and $S_{u_{l+1}}:=N_{\h' \setminus \F}(u_{l+1})$. 
Both of the sets $S_{u_1}$ and $S_{u_{l+1}}$ are subsets of $U$. 

For a set $S$, $S\subset U$  let us denote two sets $S^-$ and $S^{--}$. 
Let us denote  $S^-:=\{u_i|u_{i+1}\in S\}$ and naturally $S^{--}:=(S^-)^-$.
We may assume  $u_{l+1} \notin S_{u_1}$ and $u_{1} \notin S_{u_{l+1}}$, otherwise we have a Berge-cycle of length $l$.  

If $S^{--}_{u_1} \cap S_{u_{l+1}}\not= \emptyset$ then we  have a Berge-cycle of length $l$. 
If we assume $u_i \in S^{--}_{u_1} \cap S_{u_{l+1}}$, a hyperedge $h_{u_1} \in \h' \setminus \F $ is a hyperedge incident with $u_1$ and $u_{i+2}$ and a hyperedge  $h_{u_{l+1}} \in \h' \setminus \F $ is a hyperedge incident with $u_{l+1}$ and $u_{i}$. 
Hyperedges $h_{u_{1}}$ and $h_{u_{l+1}}$ are distinct since we have $u_{1} \notin S_{u_{l+1}}$. 
Then the following is a Berge-cycle of length $l$, 
\[u_1,f_1,u_2,\dots,u_{i-1}, h_{i-1},u_i,h_{u_{l+1}},u_{l+1},h_l,u_l, \dots h_{i+2}, u_{i+2}, h_{u_1}, u_1.\] 
Therefore we may assume that $S^{--}_{u_1} \cap S_{u_{l+1}}= \emptyset$.
Similarly we may assume that $S^{-}_{u_1} \cap S_{u_{l+1}}= \emptyset$.


 If there is a vertex $u_j$ such that $u_{j} \notin S_{u_1}$ but $u_{j+1} \in S_{u_1}$ then $u_{j-1}\notin S^{-}_{u_1}$ and $u_{j-1}\in S^{--}_{u_1}$. Hence we have $\abs{S^{-}_{u_1} \cup S^{--}_{u_1}}\geq \abs{S^{-}_{u_1}}+1=\abs{S_{u_1}}+1\geq \floor{\frac{k-1}{2}}+1$. 
 From here we split proof of this claim in some cases, depending on the parity of $k$.

\noindent \textbf{Case 1.} Let $k$ be odd. We consider two sub-cases.

\noindent \textbf{Case 1.1.}
There exists a vertex $u_j$ such that $u_{j} \notin S_{u_1}$ and $u_{j+1} \in S_{u_1}$. Since $S_{u_1}$ and $S_{u_{l+1}}$ are  subsets of $\{u_2,u_3,\dots u_l\}$, we have $\big(S_{u_1}\cup S_{u_{l+1}} \cup S^{-}_{u_1}\cup S^{--}_{u_1}\big)\subset \{u_1,u_2,u_3,\dots u_l\}$. Finally we have, 
$k-1\geq l = \abs{\{u_1,u_2,u_3,\dots u_l\}}\geq \abs{S^{-}_{u_1}\cup S^{--}_{u_1}}+\abs{S_{u_{l+1}}}\geq  \floor{\frac{k-1}{2}}+1+\floor{\frac{k-1}{2}}=k,$
 a contradiction.
 
\noindent \textbf{Case 1.2.} 
There is no vertex $u_j$ such that $u_{j} \notin S_{u_1}$ and $u_{j+1} \in S_{u_1}$. Hence we have $\{u_2,u_3,\dots,u_{\frac{k-1}{2}}, u_{\frac{k-1}{2}+1}\} \subseteq S_{u_1}$,  since $\abs{S_{u_1}}\geq \floor{\frac{k-1}{2}}$. 
From the symmetry of the Berge-path, we also have  $\{u_{l},u_{l-1},\dots,u_{\frac{k-1}{2}+1}\} \subseteq S_{u_{l+1}}$, otherwise we are done from Case 1.1.

Since $S^{-}_{u_1} \cap S_{u_{l+1}}= \emptyset$ it follows that  $l=k-1$ and $S_{u_1}=\{u_2,u_3,\dots,u_{\frac{k-1}{2}}, u_{\frac{k-1}{2}+1}\}$ and $S_{u_{l+1}}=\{u_{l},u_{l-1},\dots,u_{\frac{k-1}{2}+1}\}$. 
For every $u_i$, $u_i\in \{u_2,u_3,\dots,u_{\frac{k-1}{2}}\}$, we have an $h_{u_1}$, $h_{u_1}\in E(\h')\setminus \F$, incident with $u_1$ and $u_{i+1}$. So we have a longest Berge-path
\[u_i,f_{i-1},u_{i-1}, \dots ,f_{1},u_1,h_{u_1}, u_{i+1}, f_{i+1}, \dots, f_l,u_{l+1}.\]
Therefore either we have $N_{\h'}(\{u_1,u_2, \dots,u_{\frac{k-1}{2}}\})=\{u_1,u_2, \dots,u_{\frac{k-1}{2}},u_{\frac{k-1}{2}+1}\}$ or we have a Berge-cycle of length $l+1$ and we are done. 
Hence the vertex $u_{\frac{k+1}{2}}$ is a cut vertex of $\h'$ and all $r$-subsets of the set $\{u_1,u_2, \dots,u_{\frac{k-1}{2}},u_{\frac{k+1}{2}}\}$ are hyperedges in $\h'$, from the minimum degree condition.
Finally if $\h'$ contains a Berge-cycle of length at least $\frac{k+3}{2}$ then it's hyperedges are disjoint from the set $\{u_1,u_2, \dots,u_{\frac{k-1}{2}}\}$.
From the connectivity of $\h'$, we have a Berge-path containing all, at least $\frac{k+3}{2}$ vertices of this long Berge-cycle and $\{u_1,u_2, \dots,u_{\frac{k-1}{2}}\}$, a contradiction. If $\h'$ is Berge-cycle free of length at least $\frac{k+3}{2}$, then  from  Theorem \ref{t1} for $\h'$, we have a contradiction-
\[e(h)\leq e(\h')+\left({\binom{\frac{k}{2}-1}{r-1}}-1\right)(n-v(\h'))< e(\h_{n,k}). \]
Hence if $k$ is odd, $\h'$ contains a cycle of length $l$ or $l+1$.

\noindent \textbf{Case 2.} Let $k$ be even. 
We have  $\abs{S_{u_1}}\geq \frac{k-2}{2}$,  $\abs{S_{u_{l+1}}}\geq \frac{k-2}{2}$ and $S_{u_1},S_{u_{l+1}}\subset U$.

Assume there exist two distinct vertices $u_x$ and $u_y$ such that $u_{x},u_y \notin S_{u_1}$ and $u_{x+1}, u_{y+1} \in S_{u_1}$,  then $u_{x-1}\notin S^{-}_{u_1}$ and  $u_{y-1}\notin S^{-}_{u_1}$, but  $u_{x-1}\in S^{--}_{u_1}$ and  $u_{y-1}\in S^{--}_{u_1}$. 
Then we have $\abs{S^{-}_{u_1}\cup S^{--}_{u_1}} \geq  \frac{k-2}{2}+2$.
Since  $S_{u_1}, S_{u_{l+1}} , S^{-}_{u_1}, S^{--}_{u_1}\subset \{u_1,u_2,u_3,\dots,u_l\}$ and $ S^{-}_{u_1} \cup S^{--}_{u_1}$ is disjoint with $S_{u_{l+1}}$, we have a contradiction
\[k-1\geq l = \abs{\{u_1,u_2,u_3,\dots,u_l\}}\geq \abs{S^{-}_{u_1}\cup S^{--}_{u_1}}+\abs{S_{u_{l+1}}}\geq  \frac{k-2}{2}+2+\frac{k-2}{2}=k.\]
Therefore there are no two such vertices. To finish the proof of this case, we are going to consider all possibilities for the set $S_{u_1}$ and $S_{u_{l+1}}$.

\noindent \textbf{Case 2.1.} 
If  $S_{u_1}=\{u_2,u_3,\dots,u_{\frac{k}{2}}\}$ then $S_{u_{l+1}}\subset \{ u_{\frac{k}{2}},  u_{\frac{k+2}{2}}, \dots ,u_l\}$, since $S^{-}_{u_1} \cap S_{u_{l+1}}= \emptyset$. 
We have that $\abs{S_{u_{l+1}}}\geq \frac{k-2}{2}$ and $l<k$. Hence either $S_{u_{l+1}}= \{ u_{\frac{k}{2}},  u_{\frac{k+2}{2}}, \dots,,u_l\}$ or  $k=\ell+1$ and $S_{u_{l+1}}= \{ u_{\frac{k}{2}},  u_{\frac{k+2}{2}}, \dots ,u_{k-1}\}\setminus\{u_x\}$ for an integer $x$ satisfying $\frac{k}{2}\leq x \leq \ell $.

If $S_{u_{k}}= \{ u_{\frac{k}{2}},  u_{\frac{k+2}{2}}, \dots ,u_{l}\}$, from a similar argument as in  Case 1.2 we will get a contradiction.

If $S_{u_{l+1}}= \{ u_{\frac{k}{2}+1},  u_{\frac{k}{2}+2}, \dots ,u_k\}$, then every vertex $v$, $v\in U\setminus \{u_{\frac{k}{2}}, u_{\frac{k}{2}+1}\}$, can be a terminal vertex of a longest Berge-path with the same defining  vertex set.  We have a contradictions from a similar argument as in Case 1.2.

If $S_{u_{k}}= \{ u_{\frac{k}{2}},  u_{\frac{k}{2}+1}, \dots ,u_{k-1}\}\setminus \{u_x\}$ and 
$x\in \{\frac{k}{2}+1,\dots, k-1\}$, then every vertex $v$,  $v\in \{ u_{\frac{k}{2}+1},  
u_{\frac{k}{2}+2}, \dots, u_{k}\}$ can be a terminal vertex of a longest Berge-path. 
Since for all $i$, $i\in \{\frac{k+2}{2},\dots,k-1\}$, $i\not=x+1$, we have a hyperedge $h_{u_{k}}$ incident with $u_{k}$ and $u_{i-1}$, therefore the following is a longest Berge-path with desired property  
\[u_1,f_1,u_2,f_2, \dots,u_{i-1}, h_{u_{k}}, u_{k},f_{k-1},u_{k-1},\dots,u_i.\]
For $i=x+1$, we have $\displaystyle d(u_x)\geq  {\binom{\frac{k}{2}-1}{r-1}}$ and $u_x$ is a terminal vertex of a longest Berge-path. Hence $N_{\h'}(u_x)\subset  \{ u_{\frac{k}{2}},  u_{\frac{k}{2}+1}, \dots ,u_l\}\setminus \{u_x\} $. There exists a Berge-path  $u_j,h_{u_j},u_x,h_{u_{j+1}},u_{j+1}$ in $\h'\setminus \F$, for some $j\geq \frac{k}{2}$.
We have a longest Berge-path with $u_{x+1}$ as a terminal vertex 
\[u_1,f_1,u_2, \dots,u_{x-1}, h_{u_{x-1}}, u_{k},f_{k-1},u_{k-1},\dots,u_{x+1}\]
where the vertex $u_x$ is in between  $u_j$  and $u_{j+1}$ in the following way, $u_j,h_{u_j},u_x,h_{u_{j+1} },u_{j+1}$. 
Finally we have, all vertices $\{u_{\frac{k}{2}+1}, \dots ,u_{k-1},u_{k}\}$ as a terminal vertex. Hence either we have $N_{\h'}(\{u_{\frac{k}{2}+1}, \dots ,u_{k-1},u_{k}\})=\{ u_{\frac{k}{2}},  u_{\frac{k}{2}+1}, \dots ,u_{k}\}$ and a similar argument as in Case 1.2 will get a contradiction or we have a Berge-cycle of length at least $l$.

\noindent \textbf{Case 2.2.} 
If there exists an integer $y$ and $x$, such that  $S_{u_1}=\{u_2,u_3,\dots, u_{y}\}\setminus \{u_x\}$  then $S_{u_1}=\{u_2,u_3,\dots, u_{\frac{k+2}{2}}\}\setminus {u_x}$ and $S_{u_{l+1}}\subset \{ u_{\frac{k+2}{2}}, \dots ,u_{k}\}$, since $S^{-}_{u_1} \cap S_{u_{l+1}}= \emptyset$. 
From the symmetry of a Berge-path, we have a contradiction since this case is the same as Case 2.1.

\noindent \textbf{Case 2.3.}
  $S_{u_1}=\{u_2,u_3,\dots, u_{x}, u_{z+2},u_{z+3},\dots,{u_y}\}$,  $S_{u_{l+1}}=\{u_x,u_{x+1},\dots, u_{z-1}, u_{y},u_{y+1},\dots,$ ${u_{l}}\}$ and $l=k-1$. 
  In this case, all vertices are terminal vertices of a longest Berge-path (will be explained in the next paragraph), hence their neighborhood must be a subset of $U$. We have a contradiction since $\h'$ is connected and $v(\h')\geq r(k-1)$.
  
  By a similar argument as before, it is easy to see every vertex can be a terminal vertex for a longest Berge-path from $U\setminus\{u_x,u_y\}$. 
  We have distinct non-defining hyperedges $h_{u_1}$ incident with $u_1$ and $u_y$ and $h_{u_{l+1}}$ incident with $u_{l+1}$ and $u_{x+1}$, the following is a longest Berge-path starting at $u_x$.
  \[u_x,f_{x-1},u_{x-1},\dots,u_1, h_{u_1},u_y, f_{y-1},u_{y-1}, \dots, u_{x+1},h_{u_{l+1},f_l,u_l,\dots, f_y,u_{y+1}}.\]
  Similarly for the vertex $u_y$ we have a longest Berge-path starting at $u_y$.
  
 Finally we are done investigating all possible cases. Hence we have a Berge-cycle of length $l$ or $l+1$.
\end{proof}

\begin{claim}\label{cycle_l+_free}
 $\h'$ does not contain a Berge-cycle of length $l+1$, in other words $\h'$ is $BC_{>l}$-free.
\end{claim} 
\begin{proof}[Proof of Claim \ref{cycle_l+_free}]
Assume by contradiction that we have a Berge-cycle of length $l+1$,  since $\h'$ is connected  and $r(l+1)-(l+1)\leq rk-k<rk-r \leq v(\h')$ there exists a  vertex in $\h'$ not incident to any of the defining hyperedges of this Berge-cycle. 
Therefore there exists a  Berge-path containing  this vertex and all defining vertices of the Berge-Cycle of length $l+1$ as defining vertices. 
We got a contradiction.
\end{proof}

From the Claim \ref{Contains_cycle} and Claim \ref{cycle_l+_free} we have that  $\h'$ contains a Berge-cycle of length $l$. At this point we are ready to prove that $\h'=\h_{v(\h'),k}$.

\begin{claim}\label{cycle_l_free}
We have  $\h'=\h_{v(\h'),k}$.
\end{claim} 
\begin{proof}[Proof of Claim \ref{cycle_l_free}]
Let $\C_l$ be a Berge-cycle of length $l$  and let the vertices and hyperedges of the Berge-cycle $\C_l$ be $V(\C_l):=\{v_1,v_2, \dots, v_l\}$ and $E(\C_l):=\{h_1, h_2, \dots, h_l\}$ respectively.
Since $v(\h') \geq r\cdot (k-1)$ and  $\abs{\cup^{l}_{i=1} h_i}<r\cdot (k-1)-3$, we have distinct vertices $w_1$, $w_2$ and $w_3$ not incident with any hyperedge from $E(\C_l)$. 
We have $N_{\h'}(w_1) \subset V(\C_l)$, otherwise if $N_{\h'}(w_1) \not\subset V(\C_l)$ then we have a hyperedge $h$, $h \in E(\h') \setminus E(\C_l)$ and  a vertex $w_1$, $w_1' \in V(\h') \setminus V(\C_l)$ such that $\{w_1, w_1'\}\in h$. Then from the connectivity of $\h'$ and the minimum degree condition of $w_1$ we  have a Berge-path of length at least $l+1$, containing all defining vertices of the cycle, $w_1$ and $w_1'$, a contradiction. So we have   $N_{\h'}(w_i) \subset V(\C_l)$ for all $i$, $i\in \{1,2,3\}$.

We have $\displaystyle d_{\h'}(w_1)\geq {\binom{\floor{\frac{k-1}{2}}}{r-1}}$, hence $\abs{N_{\h'}(w_1)}=\abs{N_{\h \setminus E(\C_l)}(w_1)} \geq \floor{\frac{k-1}{2}}$. 
If $r>3$ then we have, 
\begin{equation}
\displaystyle d_{\h'}(w_1)  \geq {\binom{\floor{\frac{k-1}{2}}}{r-1}} > {\binom{\floor{\frac{k-1}{2}}-1}{r-1}} +\left(\ceil{\frac{k-1}{2}}+1\right).
\end{equation}
Hence there are at least $\floor{\frac{k-1}{2}}$ vertices  sharing at least two hyperedges with $w_1$. Let us denote this set by $S_{w_1}$ when $r>3$. If $r=3$ then since we have 
\begin{equation}
\displaystyle d_{\h'}(w_i)  > {\binom{\floor{\frac{k-1}{2}}}{r-1}}> {\binom{\floor{\frac{k-1}{2}}-1}{r-1}} +\frac{\left(\ceil{\frac{k-1}{2}}+1\right)}{2},
\end{equation}
there is a set of vertices incident to $w_1$ denoted by  $S_{w_1}$, of size  $\floor{\frac{k-1}{2}}$, such that for each pair of vertices from $S_{w_1}$ there are at least two hyperedges in $\h'$ incident to $w_1$ and one of those vertices. Similarly there is a set of vertices $S_{w_2}$, for the vertex $w_2$.
Both sets $S_{w_1}$ and $S_{w_2}$ are subset of $V(\C_l)$ and have size at least $\floor{\frac{k-1}{2}}$.
If $v_t \in S_{w_1}$ then   $v_{t+1} \not\in S_{w_1}$ (where the indices are taken modulo $l$), since we have no Berge-cycle of length $l+1$.
If $v_t \in S_{w_1}$ then   $v_{t+1} \not\in S_{w_2}$ since otherwise we have the Berge-path  $w_1, h_{w_1} v_t, h_{t-1}, v_{t-1}, \dots, h_{t+1}, v_{t+1}, h_{w_2}, w_2$ of length $l+1$, a contradiction ($h_{w_1}$ and $h_{w_2}$  are distinct hyperedges incident with $w_1,v_t$ and $v_{t+1},w_2$ respectively, they can be chosen greedily).
Hence we can assume, without loss of generality, that   $S:=S_{w_1}=S_{w_2}=\{v_1, v_3, \dots, v_{2\floor{\frac{k-1}{2}}-1}\}$. 
Therefore from the minimum degree condition, for every vertex $w_i$, $i \in \{1,2,3\}$, we have all hyperedges containing $w_i$ and any $r-1$ subset of $S$ in $\h'$. 

The rest of the proof of Claim \ref{cycle_l_free} is split into the two parts, depending on the parity of $k$.

\noindent\textbf{Case 1.}
Let $k$ be odd, then $l=k-1$ since $\frac{k-1}{2}\cdot 2 \leq l<k$. Consider the Berge-cycle $C_l$, we can  exchange a vertex $v_i$, $v_i \in \{v_2, v_4, \dots, v_{k-2}\}$,  with a vertex $w_j$,  $j \in \{1,2,3\}$. Obviously,  hyperedges $h_i$, $h_{i-1}$ will be changed automatically. Their substitutes exist since there exists every hyperedge containing $w_j$ and any $r-1$ set from $S$ and $\{v_{i-1},v_{i+1}\}\in S$.  Therefore there exists a  Berge-cycle of length $l$ not using a vertex $v_i$ and hyperedges $\{h_{i-1},h_i\}$ as defining vertices/hyperedges, for each  $i \in \{2, 4, \dots,k-2\}$. 

For every hyperedge $h$, $h \in E(\h')$, we have that $\abs{h \cap S} \geq r-1$. 
Assume a contradiction, there are $x_1$ and $x_2$ vertices, $\{x_1,x_2\} \subset h \setminus S$, for some $ h \in E(\h')$.
We may assume, without loss of generality, that $h \in E(\h')\setminus E(\C_l)$ and $x_1,x_2\notin V(\C_l)$ from the switching argument in the previous paragraph. 
If after removing $h$ from $\h'$, $x_1$ or $x_2$ are in the same connected component as $\C_l$ then we have a contradiction since we have a Berge-path longer than $l$, which contains all vertices of  $V(\C_l)$, $x_1$ and $x_2$ (the hyperedge $h$ will connect vertices $x_1$ with $x_2$).
But if $x_1$ and $x_2$  are in a different connected component(s) than  $\C_l$ after removing hyperedge $h$, then we have a contradiction since we have a Berge-path which uses all vertices from $V(\C_l)$, $x_1$ and some other vertex in its   connected component, hence a Berge-path longer than $l$. 
Therefore we have that  $\h'$ is sub-hypergraph of $\h_{v(\h'),k}$, considering set $S$ as the set $A$ in the construction of $\h_{v(\h'),k}$. On the other hand we have $e(\h')\geq e(\h_{v(\h'),k})$, hence $\h'=\h_{v(\h'),k}$.

\noindent \textbf{Case 2.} Let  $k$ be even, then $\frac{k-2}{2}\cdot 2 \leq l<k$, hence $l=k-1$  or $k-2$. If $l=k-2$, then as in Case 1, we will get $\h'=\h_{v(\h'), k-1}$ which leads us to a contradiction since $e(\h')\geq e(\h_{v(\h'),k})>e(\h_{v(\h'),k-1})$. 
Hence, we have that $l=k-1$, $S=\{v_1, v_3, \dots, v_{k-3}\}$. 

We have that for every hyperedge  $h \in E(\h')$,  either $\abs{h \cap S} \geq r-1$ or $\abs{h \cap S} = r-2$ and $h\setminus S= \{v_{k-2},v_{k-1}\}$. 
Assume by contradiction there is a hyperedge $h$ and suppose $x_1$ and $x_2\notin S$, $\{x_1,x_2\} \subset h \in E(\h')$, for some vertices $x_1,x_2$ such that  $x_1 \notin \{v_{k-2},v_{k-1}\}$.
We may assume without loss of generality that $x_1\not=w_3$.
If $h \in \{h_{k-1},h_1, h_2, \dots, h_{k-3}\}$, say $h=h_i$, then $v_{i+1}$ or $v_i$ is from $S$, without loss of generality assume $v_{i+1}$ is from $S$. 
We have a contradiction since a Berge-path  
$x_1, h_i, v_i, h_{i-1}, v_{i-1} \dots, h_{i+1}, v_{i+1}, h_{w_3}, w_3$ is longer than $l=k-1$
(where $h_{w_3}$ is a hyperedge incident with $v_{i+1}$ and $w_3$).

If the hyperedge $h$ is not a defining hyperedge and $x_1$ is a defining vertex of $\C_l$, assume $x_1=v_i$, where $v_i\in \{v_2, v_4, \dots, v_{k-2}\}$, then there exists a  Berge-cycle of length $k-1$ not using a vertex $v_i$. 
Simply we can exchange $v_i$ with $w_1$  in the Berge-cycle $C_l$.
Hence we may assume $x_1, x_2 \notin \{v_1,v_2, \dots,v_{k-3},w_3\}$.
 If $x_2\in \{v_{k-2},v_{k-1}\}$ then without loss of generality we may assume $x_2=v_{k-1}$ then we have a contradiction since there is the Berge-path 
$x_2, h, v_{k-1},h_{k-2},v_{k-2}, \dots ,h_1,v_1, h_{w_3}, w_3$ of length $l$.
Therefore from here, we may assume both $x_1, x_2$ are non-defining vertices. If the hypergraph $\h'\setminus \{h\}$ is not connected, then there is a Berge-path of length at least $l$ in $\h$, using all defining vertices of the Berge-cycle $\C_l$ and at least two vertices from another connected component, a contradiction. 
If   $\h'\setminus \{h\}$ is connected, there exists a Berge-path of length $l$ in $\h'\setminus \{h\}$, using all defining vertices of Berge-cycle $\C_l$, $x_1$ and $x_2$ as defining vertices, a contradiction. 
Finally, if $h=h_{k-2}$, and there is another hyperedge $h'$, $\{v_{k-2},v_{k-1}\}\in h'$ we could exchange $h$ with $h'$ and $h$ will no longer be a defining hyperedge so the previous argument would lead us to a contradiction.
If there is no other hyperedge $h'$, $\{v_{k-2},v_{k-1}\}\in h'$ then we have a contradiction 
\[e(\h_{v(\h'),k-1})+1\geq e(\h')\geq e(\h_{v(\h'),k})>e(\h_{v(\h'),k-1})+1\]
The first inequality holds because of the structure of $\h'$, considering set $S$ as the set $A$ in the construction. 
The second inequality holds because of the definition of $\h'$, and the third inequality holds from  the definition of the  hypergraphs $\h_{n,k}$. We have shown that for every hyperedge $h$, $h \in E(\h')$, we have $\abs{h \cap S} \geq r-1$ or $\abs{h \cap S} \geq r-2$ and $h\setminus S= \{v_{k-2},v_{k-1}\}$. 
Hence we have $\h' \subset \h_{v(\h'),k}$, and as $e(\h')\geq e(\h_{v(\h'),k})$, we have $\h'=\h_{v(\h'),k}$. 
\end{proof}

From the last Claim we may assume $\h'=\h_{v(\h'),k}$.  The following claim will finish the proof.

\begin{claim}\label{h=h'}
 We have   $\h=\h_{v(\h'),k}$.
\end{claim} 
\begin{proof}[Proof of Claim \ref{h=h'}]
We are going to show that $\h=\h'$. We got $\h'$ from $\h$ by deleting vertices from $\h$ with degree at most ${\binom{\floor{\frac{k-1}{2}}}{r-1}}-1$. We have
\begin{displaymath}
e(\h)-e(\h')\leq \bigg({\binom{\floor{\frac{k-1}{2}}}{r-1}}-1\bigg)\cdot (v(\h)-v(\h')).
\end{displaymath}
We have that $\h'=\h_{v(\h'),k}$, then
\begin{displaymath}
e(\h')={\binom{\floor{\frac{k-1}{2}}}{r-1}}\bigg(v(\h')-\floor{\frac{k-1}{2}}\bigg)+{\binom{\floor{ \frac{k-1}{2}}}{r}}+ \mathbb{1}_{2|k}{\binom{\floor{\frac{k-1}{2}}}{r-2}}.
\end{displaymath}
From the assumption we have
\begin{displaymath}
e(\h)\geq {\binom{\floor{\frac{k-1}{2}}}{r-1}}\bigg(v(\h)-\floor{\frac{k-1}{2}}\bigg)+{\binom{\floor{ \frac{k-1}{2}}}{r}}+ \mathbb{1}_{2|k}{\binom{\floor{\frac{k-1}{2}}}{r-2}}.
\end{displaymath}
Finally we have
\begin{displaymath}
 {\binom{\floor{\frac{k-1}{2}}}{r-1}}\bigg(v(\h)-v(\h')\bigg)
\leq e(\h)-e(\h') \leq
 \bigg({\binom{\floor{\frac{k-1}{2}}}{r-1}}-1\bigg)\cdot (v(\h)-v(\h')). 
\end{displaymath}
Hence we have $v(\h)=v(\h')$. Hence we got the desired result $\h=\h'$.
\end{proof}

From the final claim we have $\h=\h_{v(\h'),k}$, therefore the proof of Theorem \ref{main_theorem} is finished.

\end{proof}

\section{Remarks and Open Questions}

It would be of interest to further investigate $ex^{conn}_r(n,\B P_k)$, for $k < 2r+13$. In particular, we are interested to know what is the minimum $k$, for which $H_{n,k}$ is the only $n$-vertex, extremal, connected, $\B P_k$-free $r-$graph for large enough $n$. Moreover one may ask the stability version of Theorem \ref{main_theorem} i.e, to find the minimum number of hyperedges such that every  connected, $n$-vertex $r$-graph without a Berge-path of length $k$, with that many hyperedges is a sub-hypergraph of $H_{n,k}$. 

The other possible ways to continue this research is to ask the similar questions for the Berge-trees see Definition 1, in \cite{GSTZ} instead of paths. As one may expect those questions are harder and widely open. The related results can be found  in~\cite{bigk,GSTZ}.

\section*{Acknowledgment}

We thank the anonymous referees for their useful suggestions. The research of authors was partially supported by the National Research, Development and Innovation Office NKFIH, grants K132696, K116769, K117879 and K126853. The research of the second author is partially supported by  the Shota Rustaveli National Science Foundation of Georgia SRNSFG, grant number DI-18-118.

\end{document}